\documentclass[12pt,reqno]{amsart}
\usepackage{amssymb,amsmath,amsthm,newlfont,enumerate,color,epsfig,graphicx}
\newtheorem{prop}{Proposition}
\newtheorem*{theorem*}{Theorem}
\newtheorem*{remark*}{Remark}

\newtheorem*{definition*}{Definition}
\newtheorem{theorem}[prop]{Theorem}
\newtheorem{proposition}[prop]{Proposition}
\newtheorem{corollary}[prop]{Corollary}

\theoremstyle{remark}

\newcommand{\CC}{{\mathbb C}}
\newcommand{\C}{{\mathbb C}}

\newcommand{\cF}{{\mathcal{F}_{\alpha}^2}}
\newcommand{\dF}{{\mathcal{F}_{\alpha}^\infty}}

\newcommand{\Hol}{\operatorname{Hol}}

\begin{document}

\title[Multiple sampling and  interpolation]{Multiple sampling and  interpolation in the classical Fock space}

\author{A. Borichev, A. Hartmann, K. Kellay, X.  Massaneda }
\address{A.Borichev: Aix Marseille Universit\'e\\
CNRS\\ Centrale Marseille\\ I2M\\ 13453 Marseille\\ France \&\ Saint Petersburg State University, Department of Mathematics and Mechanics, 28 Universitetskii prosp., Staryi Petergof, 198504, Russia}
\email{alexander.borichev@math.cnrs.fr}

\address{A.Hartmann:   Universit\'e de Bordeaux\\
 IMB\\ 351 cours de la Lib\'era\-tion\\ 33405 Talence\\ France}
\email{Andreas.Hartmann@math.u-bordeaux1.fr}

\address{K.Kellay:   Universit\'e de Bordeaux\\
 IMB\\ 351 cours de la Lib\'eration\\ 33405 Talence\\ France}
\email{karim.kellay@math.u-bordeaux1.fr}

\address{X.Massaneda: Universitat  de Barcelona\\
Departament de Ma\-te\-m\`a\-tica Aplicada i An\`alisi\\
Gran Via 585, 08007-Bar\-ce\-lo\-na\\ Spain}
\email{xavier.massaneda@ub.edu}

\begin{abstract}
We study multiple sampling, interpolation and uni\-queness for the classical Fock space in the case of unbounded multiplicities. 
\end{abstract}

%\date{22/05/15}
\thanks{The work was supported by Russian Science Foundation grant 14-41-00010.}

\maketitle

Sampling and interpolating sequences in Fock spaces were characterized by Seip and Seip--Wallst\'en in \cite{S1,S2} by means of a certain Beurling--type asymptotic uniform density.   
The case of uniformly bounded multiplicities was considered by Brekke and Seip \cite{BS} who gave a complete description in this situation. Their conditions show that it is not possible that a sequence is simultaneously sampling and interpolating. 

Brekke and Seip also asked whether there exist sequences which are simultaneously sampling and interpolating when the multiplicities are 
unbounded. 

In this research note we formulate some conditions (of geometric nature) for sampling and interpolation. They show that the answer to this question is negative when the multiplicities tend to infinity. 

The detailed version of this work will be published elsewhere.
\\

We now introduce the necessary notation. For $\alpha>0$, define the Fock space  $\cF$ by 
\[
 \cF=\biggl\{f\in \Hol(\C): \|f\|^2_2=\|f\|_{\alpha,2}^2:=\frac{\alpha}{\pi}\int_\C|f(z)|^2e^{-\alpha|z|^2}dm(z)<\infty\biggr\}.
\]
The  space $\cF$ is a Hilbert space with the inner product 
$$
\langle f, g\rangle =\frac{\alpha}{\pi}\int_\CC f(z)\overline{g(z)}e^{-\alpha|z|^2}dm(z).
$$
The sequence 
\[e_k  (z)=\frac{\sqrt{\alpha^k}}{\sqrt{k!}}z^k, \qquad k\geq 0,
\] defines an orthonormal basis in $\cF$.

Recall that the translations
\[
T_z f(\zeta) =T^\alpha_zf(\zeta):=e^{\alpha \bar z\zeta-\frac{\alpha}{2}|z|^2}f(\zeta-z),\qquad f\in \cF.
\]
act isometrically in $\cF$. 
\\

Let us now define sampling and interpolation in the unbounded multiplicity case.
Consider the divisor $X=\{(\lambda,m_\lambda)\}_{\lambda\in\Lambda}$, where $\Lambda$ is a sequence of points in $\C$ and $m_\lambda\in\mathbb N$ is the
multiplicity associated with $\lambda$. 
\vspace{1em}
The divisor $X$ is called 
\begin{itemize}
\item \emph{sampling} for $\cF$ if
\[
 \|f\|^2_{2}\asymp 
  \sum_{\lambda\in \Lambda}\sum_{k=0}^{m_\lambda-1}|\langle f,T_{\lambda} e_k  \rangle|^2 ,\qquad f\in \cF\ ,
\]
\item  \emph{interpolating} for $\cF$ if for every sequence 
\[
v=\{v_\lambda^{(k)}\}_{\lambda\in\Lambda,\; 0\leq k<m_\lambda}
\] such that
\begin{equation*}
\|v\|_2^2:=\sum_{\lambda\in \Lambda}\sum_{k=0}^{m_\lambda-1}|v_\lambda^{(k)}|^2<\infty,
\end{equation*}
there exists a function $f\in \cF$  satisfying 
\begin{equation*}\label{interpol}
\langle f,T_{\lambda} e_k  \rangle =v_\lambda^{(k)},\qquad 0\leq k <  m_\lambda,\quad \lambda\in \Lambda. \\
\end{equation*}
\end{itemize}
\vspace{1em}

As in the situation of classical interpolation and sampling, the separation between points in $\Lambda$ plays an important role.
Denote by $D(z,r)$ the disc of radius $r$ centered at $z$.

\begin{itemize}
\item A divisor $X=\{(\lambda,m_\lambda)\}_{\lambda\in\Lambda}$ is said to satisfy the 
\emph{finite overlap condition} if 
\[
\sup_{z\in \CC}\sum_{\lambda\in \Lambda}\chi_{D(\lambda,\sqrt{m_\lambda/\alpha})}(z)<\infty.
\]
\end{itemize}
 
If $\Lambda$ is a finite union of $\Lambda_j$ such that the discs $D(\lambda,\sqrt{m_\lambda/\alpha})$, $\lambda\in \Lambda_j$, are disjoint for every $j$, then $X$ satisfies the finite overlap condition.\\

The following  result gives geometric conditions for sampling  in the case of unbounded multiplicities. 

\begin{theorem}\label{thm2} $\mathrm{(a)}$ If $X=\{(\lambda,m_\lambda)\}_{\lambda\in\Lambda}$ is a sampling divisor for $\cF$, then $X$ satisfies the finite overlap condition and there exists $C> 0$  such that 
\[
 \bigcup_{\lambda\in \Lambda} D(\lambda,\sqrt{m_\lambda/\alpha}+C)= \C.
\]
\\
$\mathrm{(b)}$ Conversely, if $X=\{(\lambda,m_\lambda)\}_{\lambda\in\Lambda}$ satisfies the 
finite overlap condition and for every   $C>0$ there is a  compact subset $K$ of $\C$ such that
\[
 \bigcup_{\lambda\in \Lambda, \,   m_\lambda>\alpha C^2} D(\lambda,\sqrt{m_\lambda/\alpha}-C)= \C \setminus K, 
\]
then $X$ is a sampling divisor for $\cF$.
\end{theorem}

\begin{remark*} Let $X=\{(\lambda,m_\lambda)\}_{\lambda\in\Lambda}$  satisfy the conditions of Theorem~\ref{thm2} (b). Then we can find a subset $\Lambda_1\subset \Lambda$ such that for every $C>0$ there is a compact subset $K$ of $\CC$ satisfying 
 \[
 \bigcup_{\lambda\in \Lambda_1, \,  m_\lambda>\alpha C^2} D(\lambda,\sqrt{m_\lambda/\alpha}-C)= \C \setminus K, 
\]
and 
\[
\lim_{\lambda\in \Lambda_1,\,|\lambda|\to \infty}m_\lambda=+\infty.
\]
\end{remark*}
\vspace{1em}

The following result gives geometric conditions for interpolation in the case of unbounded multiplicities.

\begin{theorem}\label{thm-int} $\mathrm{(a)}$ If  $X=\{(\lambda,m_\lambda)\}_{\lambda\in\Lambda}$ is an interpolating divisor for $\cF$, then there exists $C>0$ such that the discs $\{D(\lambda,\sqrt{m_\lambda/\alpha}-C)\}_{\lambda\in\Lambda,\,m_\lambda>\alpha C^2}$ are pairwise disjoint.
\\
$\mathrm{(b)}$ Conversely, if  the disks $\{D(\lambda,\sqrt{m_\lambda/\alpha}+C)\}_{\lambda\in\Lambda}$ are pairwise disjoint for some $C>0$, then $X$ is an interpolating divisor for $\cF$.
\end{theorem}

\vspace{1em}

\begin{remark*} It is easily seen that if $X$ is a divisor such that for some $C>0$ the discs $\{D(\lambda,\sqrt{m_\lambda/\alpha}-C)\}_{\lambda\in\Lambda, m_\lambda>\alpha C^2}$ are pairwise disjoint, and if $\lim\limits_{|\lambda|\to \infty}m_\lambda=+\infty$, then $X$ satisfies the finite overlap condition.
\end{remark*}

Though our geometric conditions do not characterize interpolation and sampling they allow us to deduce the following result, which gives a partial answer to the question raised by Brekke and Seip.
\begin{corollary}\label{thm0}
Let the divisor $X=\{(\lambda,m_\lambda)\}_{\lambda\in\Lambda}$ satisfy the condition  $\lim\limits_{|\lambda|\to \infty}m_\lambda=+\infty$.
Then $X$ cannot be simultaneously interpolating and sampling for $\cF$.
\end{corollary}

The problems of sampling and interpolation are linked to that of uniqueness, and thus to zero divisors,
for which some conditions are discussed in \cite{L, Zhu}. 

We will formulate here a necessary condition for zero divisors which, apparently, does not follow from those known so far.

\begin{proposition}\label{thm1}
Let $X=\{(\lambda,m_\lambda)\}_{\lambda\in \Lambda}$. If there exists a compact subset $K$ of $\C$ such that
\[
 \bigcup_{\lambda\in \Lambda} D(\lambda,\sqrt{m_\lambda/\alpha})=\C\setminus K,
\]
then $X$ is not a zero divisor for $\cF$. 
\end{proposition}

As a matter of fact, this result holds more generally in the weighted Fock space with uniform norm $\dF$ defined by
\[
 \dF=\left\{f\in \Hol(\CC)\text{ : } \|f\|_\infty=\|f\|_{\alpha,\infty}:=\sup_{z\in \CC}|f(z)|
 e^{-\frac{\alpha}{2}|z|^2}<\infty\right\}.
\]

\end{document}